\documentclass{article}
\usepackage{pictex,amssymb,amsmath,amsthm}

\theoremstyle{definition}
\newtheorem{thm}{Theorem}[section]
\newtheorem{prop}{Proposition}[section]
\newtheorem{lem}{Lemma}[section]
\theoremstyle{remark}

\numberwithin{equation}{section}

\newcommand{\h}{\widehat}
\newcommand{\pd}[2]{\dfrac{\partial #1}{\partial #2}}
\newcommand{\eps}{\varepsilon}
\newcommand{\de}{_{\eps}}
\newcommand{\bs}{\boldsymbol}
\newcommand{\ue}{^\eps}
\newcommand{\dx}{\,d\bs x}
\newcommand{\dhx}{\,d\widehat{\bs x}}

\newcommand{\dxi}{\,d\xi}

\newcommand{\dsigma}{\,d\sigma}
\newcommand{\LL}{{\boldsymbol L}}
\newcommand{\build}[3]{\mathrel{\mathop{\kern 0pt#1}\limits_{#2}^{#3}}}

\def\curl{\operatorname{\mathbf{curl}}}
\def\diver{\operatorname{div}}
\def\Int{\operatorname{Int}}
\newcommand{\vp}{\times}

\begin{document}

\title{Asymptotic behaviour of the inductance coefficient for thin conductors}

\author{\sc Youcef Amirat, Rachid Touzani\\
{\small\sl Laboratoire de Math\'ematiques Appliqu\'ees, UMR CNRS 6620}\\
{\small\sl Universit\'e Blaise Pascal (Clermont--Ferrand)}\\
{\small\sl 63177 Aubi\`ere cedex, France}}
\date{}

\maketitle
\begin{abstract}
We study the asymptotic behaviour of the inductance coefficient for a thin toroidal inductor
whose thickness depends on a small parameter $\eps>0$. We give an explicit form of the
singular part of the corresponding potential $u\ue$ which allows to construct the limit potential
$u$ (as $\eps\to 0$) and an approximation of the inductance coefficient $L\ue$.
We establish some estimates of the deviation $u\ue-u$ and of the error of approximation of
the inductance. We show that $L\ue$ behaves asymptotically as $\ln\eps$, when $\eps\to 0$.

\bigskip
\begin{center}{\bf R\'esum\'e}\end{center}
On \'etudie le comportement asymptotique du coefficient d'inductance pour un inducteur toro\"\i{}dal
filiforme dont l'\'epaisseur d\'epend d'un petit para\-m\`et\-re $\eps>0$. On donne une forme
explicite de la partie singuli\`ere du potentiel associ\'e $u\ue$ puis on construit le
potentiel limite $u$ (quand $\eps\to~0$) et on donne une approximation du coefficient d'inductance
$L\ue$. On \'etablit des estimations de l'\'ecart $u\ue-u$ et de l'erreur d'approximation
de l'inductance. On montre que $L\ue$ se comporte asymptotiquement comme $\,\ln\eps\,$ au
voisinage de $\eps=0$.
\end{abstract}

\vfill
{\small
}

\eject

\section{Introduction}

Electrotechnical devices often involve thick conductors in which a magnetic field can be
induced, and thin wires or coils, as inductors, connected to a power source generator. The
problem is then to derive mathematical models which take into account the simultaneous
presence of thick conductors and thin inductors. For a two--dimensional configuration
where the magnetic field has only one nonvanishing component, it was shown that the
eddy current equation has the Kirchhoff circuit equation as a limit problem, as the thickness
of the inductor tends to zero, see \cite{Touzani}.
For the three--dimensional case, eddy current models require the use of a relevant quantity
that is the self inductance of the inductor, see \cite{Bossavit}, \cite{BossVer}.
This number has to be evaluated {\em a priori}
as a part of problem data. It is the purpose of the present paper to study the asymptotic
behaviour of this number when the thickness of the inductor goes to zero.

Let us consider a toroidal domain of $\mathbb R^3$, denoted by $\Omega\de$, whose thickness
depends on a small parameter $\eps>0$. The geometry of $\Omega\de$ will be described in the
next section. We denote by $\Gamma\de$ the boundary of $\Omega\de$, by $n\de$ the outward unit
normal to $\Gamma\de$, and by $\Omega'\de$ the complementary of its closure, that is
$\Omega'\de=\mathbb R^3\setminus\overline\Omega\de$. We denote by $\Sigma$ a cut in the
domain $\Omega'\de$, that is, $\Sigma$ is a smooth orientable surface such that, for any $\eps>0$,
$\Omega'\de\setminus\Sigma$ is simply connected.

Let now $\bs h\ue$ denote the time--harmonic and complex valued
magnetic field. Neglecting the displacement currents, it follows from Maxwell's equations
that
$$
\curl\bs h\ue = 0,\ \diver\bs h\ue = 0\quad\text{in }\Omega'\de.
$$
Then, by a result in \cite{Dautray-Lions-2}, p. 265, $\bs h\ue$ may be written in the form
\begin{equation}
\bs h\ue_{|\Omega'\de} = \nabla\varphi\ue + I\ue\nabla u\ue,\label{h-ext}
\end{equation}
where $I\ue$ is a complex number, $\varphi\ue\in W^1(\Omega'\de)$ and satisfies
$$
\Delta\varphi\ue = 0\qquad\text{in }\Omega'\de,
$$
and $u\ue$ is solution of~:
\begin{equation}
\left\{\begin{aligned}{}
&\Delta u\ue = 0&&\qquad\text{in }\Omega'\de\setminus\Sigma,\\
&\pd{u\ue}n = 0&&\qquad\text{on }\Gamma\de,\\
&[u\ue]_{\Sigma} = 1,\\
&\left[\pd{u\ue}n\right]_\Sigma = 0.
\end{aligned}\right.\label{Pb-u}
\end{equation}

Here $W^1(\Omega'\de)$ is the Sobolev space
$$
W^1(\Omega'\de) = \left\{v;\ \rho v\in L^2(\Omega'\de),\ \nabla v\in \LL^2(\Omega'\de)\right\},
$$
equipped with the norm
\begin{equation}
\|v\|_{W^1(\Omega'\de)} = \left(\|\rho v \|^2_{L^2(\Omega'\de)} +
\|\nabla v\|^2_{\LL^2(\Omega'\de)}\right)^{\frac 12},\label{norm}
\end{equation}
where $\LL^p(\Omega'\de)$ denotes the space $L^p(\Omega'\de)^3$ and $\rho$ is the weight
function $\rho(\bs x) = (1+|\bs x|^2)^{-\frac 12}$. Let us note here, see \cite{Dautray-Lions-2},
pp. 649--651, that
$$
|v|_{W^1(\Omega'\de)} = \left(\int_{\Omega'\de}|\nabla v|^2\dx\right)^{\frac 12}
$$
is a norm on $W^1(\Omega'\de)$, equivalent to \eqref{norm}. In \eqref{Pb-u}, $n$ is the unit
normal on $\Sigma$, and $[u\ue]_\Sigma$ (resp. $\Big[\pd{u\ue}n\Big]_\Sigma$) denotes the
jump of $u\ue$ (resp. $\pd{u\ue}n$) across $\Sigma$.

In \eqref{h-ext}, the number $I\ue$ can be interpreted as the total current flowing in the
inductor, see \cite{BossVer}.

The inductance coefficient is then defined by the expression
\begin{equation}\label{def-L}
L\ue = \int_{\Omega'\de\setminus\Sigma}|\nabla u\ue|^2\dx.
\end{equation}

Our goal is to study the asymptotic behaviour of $u\ue$ and $L\ue$ as $\eps$ goes to zero.
We first give an explicit form of the singular part of the potential $u\ue$ which allows to
construct the limit potential $u$ (as $\eps\to 0$) and an approximation of the inductance
$L\ue$. We then prove that the deviation $\|u\ue-u\|_{W^1(\Omega'\de)}$ and the error of
approximation of $L\ue$ is at order $O(\eps^{\frac 56})$. Finally we show that the inductance
coefficient $L\ue$ behaves asymptotically as $\ln\eps$, when $\eps\to 0$, and we thus
recover the result stated (without proof) in \cite{Landau}, p. 137.

The remaining of this paper is organized as follows. In Section~2 we precise the geometry of
the inductor by considering that this one is obtained by generating a
toroidal domain around a closed curve, the internal radius of the torus being
proportional to a small positive number $\eps$. Section~3 states the main result and
Section~4 is devoted to the proof.

\section{Geometry of the domain}

We consider a toroidal domain, with a small cross section. This domain may be defined as a
tubular neighborhood of a closed curve. Let $\gamma$ denote a closed Jordan arc of class $\mathcal C^3$
in $\mathbb R^3$, with a parametric representation defined by a function $\bs g:[0,1]\to\mathbb R^3$
satisfying
\begin{equation}
\bs g(0) = \bs g(1),\ \bs g'(0) = \bs g'(1),
\ |\bs g'(s)|\ge C_0>0.\label{hyp-g}
\end{equation}
For each $s\in (0,1]$ we denote by $(\bs t(s),\bs\nu(s),\bs b(s))$
the Serret--Fr\'enet coordinates at the point $\bs g(s)$, $i.e.$, $\bs t(s),\bs\nu(s),\bs b(s)$
are respectively the unit tangent vector to $\gamma$, the principal normal and
the binormal, given by
$$
\bs t=\frac{\bs g'}{|\bs g'|},\ \bs\nu = \frac{\bs t'}{|\bs t'|},\ \bs b=\bs
t\vp\bs\nu.
$$
We have the following well-known Serret--Fr\'enet formulae~:
$$
\bs t' = \kappa\bs\nu,\ \bs\nu' = -\kappa\bs t + \tau\bs b,\ \bs b' = -\tau\bs\nu,
$$
where $\kappa$ and $\tau$ denote respectively the curvature and the torsion of the arc $\gamma$.

Let $\h\Omega=(0,1)^2\times (0,2\pi)$ and let $\delta$ denote a
positive number to be chosen in a convenient way. We define, for
any $\eps$, $0\le\eps <\delta$,  the mapping $\bs F\de:\h\Omega\to\mathbb R^3$ by
$$
\bs F\de(s,\xi,\theta) = \bs g(s) + r\de(\xi)(\cos\theta\,\bs\nu(s)+\sin\theta\,\bs b(s)),
$$
where $r\de(\xi) = (\delta-\eps)\xi + \eps$. We have
$$
\begin{aligned}{}
\pd{\bs F\de}s &= \bs g' + r\de(\cos\theta\,\bs\nu'+\sin\theta\,\bs b')\notag\\
&= (|\bs g'|-r\de\kappa\cos\theta)\bs t+r\de\tau(\cos\theta\,\bs b - \sin\theta\,\bs\nu), \notag\\
\pd{\bs F\de}\xi &= (\delta -\eps)(\cos\theta\,\bs\nu + \sin\theta\,\bs b),\notag\\
\pd{\bs F\de}\theta &= r\de(-\sin\theta\,\bs\nu + \cos\theta\,\bs b).\notag
\end{aligned}
$$

The jacobian of $\bs F\de$ is therefore given by
$$
J\de(s,\xi,\theta) = (\delta-\eps)a\de(s,\xi,\theta) r\de(\xi),
$$
where
$$
a\de(s,\xi,\theta) = |\bs g'(s)|-r\de(\xi)\kappa(s)\cos\theta.
$$

According to \eqref{hyp-g}, if $\delta$ is chosen such that
$$
\delta|\kappa(s)| < |\bs g'(s)|,\qquad 0\le s\le 1,
$$
then
\begin{equation}
0<C_1\le a\de\le C_2,\label{estim-param}
\end{equation}
and the mapping $\bs F\de$ is a $\mathcal C^1$--diffeomorphism from $\h\Omega$
into $\Lambda\de^\delta=\bs F\de(\h\Omega)$.

Here and in the sequel, the quantities $C, C_1, C_2, \dots$ denote generic positive
numbers that do not depend on $\eps$.

$$
\vbox{
\beginpicture
\normalgraphs
\setcoordinatesystem units <0.08 truecm,0.08 truecm> point at 0 0
\setplotarea x from 0 to 80, y from 0 to 80
\setlinear
\startrotation by 0.966 0.259 about 25 35
\ellipticalarc axes ratio 1:3.5   360 degrees from 25 20  center at 25 35
\ellipticalarc axes ratio 1:3.5   360 degrees from 25 10  center at 25 35
\setquadratic
\plot 0 22  10 30  25 35 /
\plot 24 60  40 66  60 67 /
\plot 26 10.3  42 16.3  62 17.3 /
\plot 26 20.3 28 21.2 31 21.5 /
\plot 24.9 50.1 27.5 50.9 30.5 51.2 /
\put {$\gamma$} at 4 30
\put {$\Omega_\varepsilon$} at 24 40
\put {$\Lambda^\delta_\varepsilon$} at 24 55
\put {$\Gamma_\delta$} at 17.5 20
\put{{\sc Figure} 1 -- A sketch of the inductor geometry} at 20 -5
\endpicture
}$$

We now set, for any $0<\eps<\delta$,
$$
\Omega_\delta =\Lambda_0^\delta = \bs F_0(\h\Omega),\ \Omega'_\delta = \mathbb R^3
\setminus{\overline\Omega}_\delta,
\ \Omega'\de=\Int(\overline\Omega'_\delta\cup\overline\Lambda\de^\delta),\
\Omega\de = \mathbb R^3\setminus\overline\Omega'\de.
$$
For technical reasons, we choose in the sequel $0<\eps\le\frac \delta 2$.

Given a function $v$ on $\Lambda\de^\delta$, we define the function $\h v$ on $\h\Omega$ by
$\h v=v\circ\bs F\de$. If $v\in L^p(\Lambda\de^\delta)$, $1\le p\le\infty$, then $\h
v\in L^p(\h\Omega)$ and we have
$$
\int_{\Lambda\de^\delta} v\,d\bs x =\int_{\h\Omega}\h v\,(\delta-\eps)\,a\de r\de\,d\h{\bs x}.
$$
If $v\in W^{1,p}(\Lambda\de^\delta)$, $1\le p\le\infty$, then $\h v\in
W^{1,p}(\h\Omega)$ and we have
\begin{align}
&\pd{\h v}s = \h{\nabla v}\cdot\pd{\bs F\de}s = \h{\nabla v}\cdot
(a\de\bs t + r\de\tau\cos\theta\,\bs b - \bs r\de\tau\sin\theta\,\bs\nu),\label{dvds}\\
&\pd{\h v}\xi = \h{\nabla v}\cdot\pd{\bs F\de}\xi = (\delta-\eps)\h{\nabla v}\cdot
(\cos\theta\,\bs\nu + \sin\theta\,\bs b),\label{dvdxi}\\
&\pd{\h v}\theta = \h{\nabla v}\cdot\pd{\bs F\de}\theta = r\de\h{\nabla v}\cdot
(-\sin\theta\,\bs\nu + \cos\theta\,\bs b).\label{dvdtheta}
\end{align}
From \eqref{dvdxi} and \eqref{dvdtheta} we deduce
\begin{align}
&\h{\nabla v}\cdot\bs b = \frac{\sin\theta}{\delta-\eps}\pd{\h v}\xi + \frac{\cos\theta}{r\de}
\pd{\h v}\theta,\label{dv-db}\\
&\h{\nabla v}\cdot\bs\nu = \frac{\cos\theta}{\delta-\eps}\pd{\h v}\xi - \frac{\sin\theta}{r\de}
\pd{\h v}\theta,\label{dv-dnu}
\end{align}
and then, with \eqref{dvds} we get
\begin{equation}
\h{\nabla v}\cdot\bs t = \frac 1{a\de}\left(\pd{\h v}s-\tau\pd{\h v}\theta\right).
\label{dv-dt}
\end{equation}
Therefore, for $u$ and $v$ in $H^1(\Lambda\de^\delta)$,
\begin{align}
\int_{\Lambda\de^\delta}\nabla u.\nabla v\,d\bs x
&= (\delta-\eps)\int_{\h\Omega}
\bigg(\dfrac{r\de}{a\de}\pd{\h u}s\pd{\h v}s +
\dfrac{r\de a\de}{(\delta-\eps)^2}\pd{\h u}\xi\pd{\h v}\xi\notag\\
&\quad\qquad\qquad\quad + \left(\dfrac{a\de}{r\de}+\dfrac{\tau^2r\de}{a\de}
\right)\pd{\h u}\theta\pd{\h v}\theta\notag\\
&\quad\qquad\qquad\quad -\dfrac{r\de\tau}{a\de}\left(\pd{\h u}s
\pd{\h v}\theta + \pd{\h u}\theta\pd{\h v}s\right)\bigg)
\,d\h{\bs x}.\label{VariableChange}
\end{align}

We also define the set $\h\Gamma=(0,1)\times(0,2\pi)$ and the mapping
$\bs G\de:\h\Gamma\to\mathbb R^3$ by
$$
\bs G\de(s,\theta) = \bs g(s) + \eps(\cos\theta\,\bs\nu(s)+\sin\theta\,\bs b(s)).
$$
The boundary of $\Omega'\de$ is then represented by $\Gamma\de = \overline{\bs G\de(\h\Gamma)}$.
We have
$$
\begin{aligned}{}
\pd{\bs G\de}s &= (|\bs g'|-\eps\kappa\cos\theta)\bs t + \eps
\tau(\cos\theta\,\bs b-\sin\theta\,\bs\nu),\\
\pd{\bs G\de}\theta &= \eps(-\sin\theta\,\bs\nu + \cos\theta\,\bs b).
\end{aligned}\notag
$$
If $w\in L^2(\Gamma\de)$, we define $\h w\in L^2(\h\Gamma)$ by $\h
w=w\circ\bs G\de$, and we have
\begin{equation}
\int_{\Gamma\de} w\,d\sigma =\int_{\h\Gamma}\h
w\,\eps(|\bs g'|-\eps\kappa\cos\theta)\,d\h\sigma.\label{Int-Gamma}
\end{equation}

Clearly, $\Omega\de$ and its complementary $\Omega'\de$ are connected domains but they
are not simply connected. To define a cut in $\Omega'\de$, we denote by $\Sigma_0$ the set
$\bs F_0((0,1)^2\!\vp\!\{0\})$ and $\partial\Sigma_0 = \bs F_0((0,1)\!\vp\!\{1\}\!\vp\!\{0\})$.
Let $\Sigma'$ denote a smooth simple surface that has $\partial\Sigma_0$ as a boundary and
such that the surface $\Sigma= \Sigma'\cup\Sigma_0$ is oriented and of class $\mathcal C^1$
(cf. \cite{Kreyszig}). We denote by $\Sigma^+$ ({\em resp.} $\Sigma^-$) the oriented
surface with positive ({\em resp.} negative) orientation, and by $n$ the unit normal
on $\Sigma$ directed from $\Sigma^+$ to $\Sigma^-$. If $w\in W^1(\mathbb R^3\setminus\Sigma)$,
we denote by $[w]_\Sigma$ the jump of $w$ across $\Sigma$ through $n$, {\em i.e.}
$$
[w]_\Sigma = w_{|\Sigma^+} - w_{|\Sigma^-}.
$$

\section{Formulation of the problem and statement of the result}

We consider the boundary value problem
\begin{equation}
\left\{\begin{aligned}{}
&\Delta u\ue = 0&&\qquad\text{in }\Omega'\de\setminus\Sigma,\\
&\pd{u\ue}{n\de} = 0&&\qquad\text{on }\Gamma\de,\\
&[u\ue]_{\Sigma} = 1,\\
&\left[\pd{u\ue}n\right]_\Sigma = 0,
\end{aligned}\right.\label{Pb-eps}
\end{equation}
where $n\de$ denotes the unit normal on $\Gamma\de$ pointing outward
$\Omega'\de$ and $\bs n$ is the unit normal on $\Sigma$ oriented from $\Sigma^+$
toward $\Sigma^-$. The inductance coefficient is defined by
\begin{equation}
L\ue = \int_{\Omega'\de\setminus\Sigma}|\nabla u\ue|^2\,d\bs x.\label{induct}
\end{equation}
We want to describe the asymptotic behaviour of $u\ue$ and $L\ue$ as $\eps\to 0$.

We first exhibit a function that has the same singularity as might have the
solution of Problem \eqref{Pb-eps} (as $\eps\to 0$). Let us define
$$
\h v(s,\xi,\theta) =
\dfrac\theta{2\pi}\,\h\varphi(\xi),\qquad(s,\xi,\theta)\in\h\Omega,
$$
where $\h\varphi\in\mathcal C^2(\mathbb R)$ and such that
$$
\h\varphi(\xi)=1\text{ for } 0\le\xi\le\frac 12,\quad
\h\varphi(\xi)=0\text{ for } \xi\ge\frac 34.
$$
We then define $v:\mathbb R^3\to\mathbb R$ by~:
$$
v(\bs x) = \begin{cases} \h v(\bs F_0^{-1}(\bs x)) &\text{if }\bs x\in\Omega_\delta,\\
0&\text{if }\bs x\in\Omega'_\delta.\end{cases}
$$
Let us also define
\begin{align}
&\h f(s,\xi,\theta) =
\frac 1{2\pi a_0}\left(\dfrac{\kappa\,\sin\theta}{\delta\xi}
- \frac{\tau^2\,\delta\,\xi\,\kappa\,\sin\theta}{a_0^2}
- \pd{}s\left(\frac \tau{a_0}\right)\right)\h\varphi,\notag\\
&\qquad\qquad\qquad\qquad\quad
 + \dfrac\theta{2\pi\, a_0\,\delta^2\,\xi}(2a_0-|\bs g'|)\,\h\varphi'
+ \frac\theta{2\pi\,\delta^2}\,\h\varphi^{\prime\prime},
\qquad (s,\xi,\theta)\in\h\Omega,\notag\\
&f(\bs x) = \begin{cases} \h
f(\bs F^{-1}_0(\bs x))&\qquad\text{if }\bs x\in\Omega_\delta,\\
0&\qquad\text{if }\bs x\in\Omega'_\delta,\end{cases}\notag\\
&\varphi(\bs x) = \begin{cases}
\h\varphi(\xi)&\qquad\qquad\quad\text{if }\bs x\in\Omega_\delta,
\text{ with }(s,\xi,\theta) = \bs F_0^{-1}(\bs x),\\
0&\qquad\qquad\quad\text{if }\bs x\in\Omega'_\delta.\end{cases}\notag
\end{align}
We have the following result.

\begin{prop}
The function $v$ is solution of
\begin{equation}
\left\{\begin{aligned}{}
&\Delta v = f&&\qquad\text{in }\mathbb R^3\setminus\Sigma,\\
&[v]_\Sigma=\varphi,\\
&\left[\pd vn\right]_\Sigma = 0.
\end{aligned}\right.\label{Singular-Pb}
\end{equation}
Moreover, it satisfies
\begin{equation}
\pd v{n\de} = 0\qquad\text{on }\Gamma\de.\label{v-on-Gamma}
\end{equation}
\end{prop}

\begin{proof}
The first equation in \eqref{Singular-Pb} follows readily from definitions of $f$ and $v$. It
remains to check the boundary conditions.
On $\Sigma'_0$, we have obviously
$$
[v]_{\Sigma'_0}=\left[\pd vn\right]_{\Sigma'_0} =0.
$$
On $\Sigma_0$, we have
$$
v_{|\Sigma^+_0} = \varphi,\ v_{|\Sigma^-_0} =0,
$$
whence $[v]_\Sigma=\varphi$.
We also have, according to \eqref{dv-db}, \eqref{dv-dnu},
$$
\begin{aligned}{}
&\h{\nabla v}_{\big|\h\Sigma_0} = -\frac{\tau}{2\pi a_0}\h\varphi\,\bs t +
\frac 1\delta\,\h\varphi'\,\bs\nu +
\frac 1{2\pi\delta\xi}\,\h\varphi\,\bs b&&\qquad\text{for }\theta=2\pi,\\
&\h{\nabla v}_{\big|\h\Sigma_0} = -\frac{\tau}{2\pi a_0}\h\varphi\,\bs t
+ \frac 1{2\pi \delta\xi}\,\h\varphi\,\bs b
&&\qquad\text{for }\theta=0,
\end{aligned}
$$
with $\h\Sigma_0=(0,1)^2$. The normal to $\Sigma_0$ is defined by
$$
\h{\bs n} = \frac 1{((|\bs g'|-\delta\xi\kappa)^2+\delta^2\xi^2\tau^2)^{\frac 12}}
\,((|\bs g'|-\delta\xi\kappa)\bs b - \delta\xi\tau\bs t).
$$
Therefore
$$
\h{\pd vn_{\big|\Sigma_0}} =
\frac{\h\varphi}{2\pi((|\bs g'|-\delta\xi\kappa)^2+\delta^2\xi^2\tau^2)^{\frac 12}}
\left(\frac{|\bs g'|-\delta\xi\kappa}{\delta\xi}+\frac{\delta\xi\tau^2}{a_0}\right),
$$
and then
$$
\left[\pd vn\right]_{\Sigma_0} = 0.
$$
We have, by \eqref{dv-db}--\eqref{dv-dt},
\begin{align}
\h{\nabla v} = \frac 1{a_0}\Big(\pd{\h v}s-\tau\pd{\h v}\theta\Big)\,\bs t
&+ \Big(\frac{\cos\theta}\delta\pd{\h v}\xi-\frac{\sin\theta}{\delta\xi}\pd{\h v}\theta\Big)
\,\bs\nu\notag\\
&+ \Big(\frac{\sin\theta}\delta\pd{\h v}\xi+\frac{\cos\theta}{\delta\xi}\pd{\h v}\theta\Big)
\,\bs b.\notag
\end{align}
The normal to $\Gamma\de$ is parametrically represented by $-(\cos\theta\bs\nu +\sin\theta\bs b)$.
Then, since $\h\varphi'(\frac\eps\delta)=0$,
$$
\h{\pd{v}{n\de}_{\big|\Gamma\de}} = -\frac 1\delta\pd{\h v}\xi(s,\frac\eps\delta,\theta)
= -\frac{\theta}{2\pi\delta}\h\varphi'(\frac\eps\delta) = 0.
$$
We conclude that $v$ is solution of Problem \eqref{Singular-Pb}.
\end{proof}

\begin{lem}\label{f}
For any $1\le p<2$ we have
$$
f\in L^p(\mathbb R^3),\ v\in L^\infty(\mathbb R^3)\cap W^{1,p}(\mathbb R^3\setminus\Sigma).
$$
\end{lem}

\begin{proof}
Clearly $v \in L^\infty(\mathbb R^3)$. Let us calculate the $L^p$--norm of $f$.
Using the mapping $\bs F_0^{-1}$, we have
\begin{align*}
\|f\|^p_{L^p(\mathbb R^3\setminus\Sigma)} &= \|f\|^p_{L^p(\Omega_\delta)}\notag\\
&= \frac 1{(2\pi)^p}\int_{\h\Omega}\bigg|\frac 1{a_0}\left(\dfrac{\kappa\,\sin\theta}{\delta\xi}
- \frac{\tau^2\,\delta\,\xi\,\kappa\,\sin\theta}{a_0^2}
- \pd{}s\left(\frac\tau{a_0}\right)\right)\h\varphi\notag\\
&\qquad\qquad\qquad\qquad
 + \dfrac\theta{a_0\,\xi\,\delta^2}(2a_0-|\bs g'|)\,\h\varphi'
+ \frac\theta{\delta^2}\,\h\varphi^{\prime\prime}\bigg|^p\,\delta^2\,a_0\,\xi\dhx.\notag
\end{align*}
Owing to \eqref{estim-param} and to the fact that $\h\varphi$ is of class $\mathcal C^2$, we deduce
that the above integral is finite provided that $1\le p<2$.

Using \eqref{dv-db}--\eqref{dv-dt}, we get
$$
\|\nabla v\|^p_{\LL^p(\mathbb R^3\setminus\Sigma)} = \frac{\delta^2}{(2\pi)^p}
\int_{\h\Omega}a_0\,\xi\left|\frac{\theta^2}{\delta^2}(\h\varphi')^2+
\left(\frac 1{\delta^2\xi^2}+\frac{\tau^2}{a_0^2}\right)\h\varphi^2\right|^{\frac p2}\dhx.
$$
With the same argument as for $f$, we deduce that the above integral is finite iff $1\le p<2$.
\end{proof}

Let us now set $w\ue=u\ue-v$. We have by subtracting \eqref{Singular-Pb} from
\eqref{Pb-eps},
\begin{equation}
\left\{\begin{aligned}{}
&-\Delta w\ue = f&&\qquad\text{in }\Omega'\de\setminus\Sigma,\\
&\pd{w\ue}{n\de} = 0&&\qquad\text{on }\Gamma\de,\\
&[w\ue]_\Sigma=1-\varphi,\\
&\left[\pd{w\ue}n\right]_\Sigma=0.
\end{aligned}\right.\label{Pb-w}
\end{equation}
We note here that
Problem \eqref{Pb-w} differs from \eqref{Pb-eps} by the value of the jump of the solution
across $\Sigma$ and by the presence of a right-hand side $f$. However, we notice that $(1-\varphi)$
vanishes in a neighborhood of $\partial\Sigma$ and then, for Problem \eqref{Pb-w}, the jump
of $w\ue$ vanishes in a neighborhood of $\partial\Sigma$.

Now, to study the asymptotic behaviour of $w\ue$ and $L\ue$ as $\eps\to 0$ we consider the
following decomposition. Let $w_1$ denote the solution of
\begin{equation}
\left\{\begin{aligned}{}
&\Delta w_1 = 0&&\qquad\text{in }\mathbb R^3\setminus\Sigma,\\
&[w_1]_\Sigma=1-\varphi,\\
&\left[\pd{w_1}n\right]_\Sigma=0,\\
&w_1(\bs x)=O(|\bs x|^{-1})&&\qquad |\bs x|\to\infty.
\end{aligned}\right.\label{Pb-w1}
\end{equation}
Using \cite{Dautray-Lions-2}, p. 654, and the fact that $(1-\varphi)$ vanishes in a neighborhood
of $\partial\Sigma$, we see that Problem \eqref{Pb-w1} has a unique solution in
$W^1(\mathbb R^3\setminus\Sigma)$ given by
\begin{equation}
w_1(\bs x) = \frac 1{4\pi}\int_\Sigma(1-\varphi(\bs y))\frac{\bs n(\bs y)\cdot(\bs
x-\bs y)} {|\bs x-\bs y|^3}\,d\sigma(\bs y),\qquad\bs x\in\mathbb R^3\setminus\Sigma.
\label{expr-w1}
\end{equation}

Then we write $w\ue=w_1+w\ue_2$, where the function $w\ue_2$ is solution
of the exterior Neumann problem~:
\begin{equation}
\left\{\begin{aligned}{}
&-\Delta w_2\ue = f&&\qquad\text{in }\Omega'\de,\\
&\pd{w_2\ue}{n\de}=-\pd{w_1}{n\de}&&\qquad\text{on }\Gamma\de,\\ &w_2\ue(\bs x)=O(|\bs
x|^{-1})&&\qquad |\bs x|\to\infty.
\end{aligned}\right.\label{Pb-w2}
\end{equation}

We have the following result.
\begin{lem}
Problem \eqref{Pb-w2} admits a unique solution $w_2\ue\in W^1(\Omega'\de)$.
\end{lem}

\begin{proof}
Differentiating \eqref{expr-w1}, we obtain for $\bs x\in\Gamma\de$~:
\begin{align}
\pd {w_1}{n\de}(\bs x) &= \frac 1{4\pi}\int_\Sigma(1-\varphi(\bs y))\frac{\bs n\de(\bs x)\cdot
\bs n(\bs y)}{|\bs x-\bs y|^3}\,d\sigma(\bs y)\notag\\
&\quad- \frac 3{4\pi}\int_\Sigma(1-\varphi(\bs y))\frac{(\bs n\de(\bs x)\cdot (\bs x-\bs y))
\,(\bs n(\bs y)\cdot (\bs x-\bs y))}{|\bs x-\bs y|^5}\,d\sigma(\bs y).\notag
\end{align}
Owing to the definition of $\varphi$, the integrals over $\Sigma$ reduce to those over
$\widetilde\Sigma$ where
$$
\widetilde\Sigma = \bs\Phi\de((0,1)\times(\frac 12,1))\cup\Sigma'.
$$
So, for $\bs x\in\Gamma\de$ and $\bs y\in\widetilde\Sigma$, $|\bs x-\bs y|\ge \frac \delta 4$
since $\eps$ is chosen not greater than $\frac \delta 2$.
Therefore
\begin{equation}
\left\|\pd {w_1}{n\de}\right\|_{L^\infty(\Gamma\de)}\le C,\label{Estim-dw1dn-inf}
\end{equation}
and, since $f_{|\Omega'\de}\in L^2(\Omega'\de)$, then Problem \eqref{Pb-w2} is a classical
exterior Neumann problem which admits a unique solution
$w_2\ue\in W^1(\Omega'\de)$, see \cite{Dautray-Lions-1}, p. 343.
\end{proof}

Let finally $w_2$ denote the unique solution in $W^1(\mathbb R^3)$ of
\begin{equation}
\left\{\begin{aligned}{}
&-\Delta w_2 = f&&\qquad\text{in }\mathbb R^3,\\
&w_2(\bs x) = O(|\bs x|^{-1}),&&\qquad |\bs x|\to\infty.
\label{eq-w2}
\end{aligned}\right.
\end{equation}

As it is classical (see \cite{Nedelec} for instance) the function $w_2$ is given by
$$
w_2(\bs x) = \frac 1{4\pi}\int_{\mathbb R^3}\frac{f(\bs y)}{|\bs x-\bs y|}\,d\bs y,\qquad
\bs x\in\mathbb R^3.
$$
Summarizing the decomposition process of the solution to Problem \eqref{Pb-eps}, we have
$$
u\ue = v + w_1 + w\ue_2\qquad\text{in }\Omega'\de\setminus\Sigma,
$$
where $v$, $w_1$ and $w\ue_2$ are solutions of \eqref{Singular-Pb}, \eqref{Pb-w1} and \eqref{Pb-w2}
respectively.

We now state our main result.

\begin{thm}\label{MainResult}
Let $u\ue$ be the solution of Problem \eqref{Pb-eps} and let $L\ue$ be the
inductance coefficient defined by \eqref{induct}. Let $u$ be the function defined in
$\mathbb R^3\setminus\Sigma$ by $u=v+w_1+w_2$, where $v$, $w_1$ and $w_2$ are
solutions of \eqref{Singular-Pb}, \eqref{Pb-w1} and \eqref{eq-w2} respectively. Then
for any $\eta>0$~:
\begin{align}
&\|u-u\ue\|_{W^1(\Omega'\de)} = O(\eps^{\frac 56-\eta}),\label{err-u}\\
&L\ue = -\dfrac{\ell_\gamma}{2\pi}\ln\eps + L' - \int_{\mathbb R^3}f(w_1+w_2)\,d\bs x\notag\\
&\qquad+ \int_\Sigma(1-\varphi)\left(\pd{w_1}n+\pd{w_2}n+2\pd vn\right)\,d\sigma +
O(\eps^{\frac 56-\eta}),\label{expr-Le}
\end{align}
where $\ell_\gamma$ is the length of the curve $\gamma$ and
$$
L' = \frac{\ell_\gamma}{2\pi}\ln{\frac\delta 2} +
\frac 1{4\pi^2}\int_{\h\Omega}\left(a_0\xi\theta^2(\h\varphi')^2 +
\frac{\delta^2\xi\tau^2}{a_0}\h\varphi^2\right)\dhx
+\ell_\gamma\int_{\frac 12}^1\frac{\h\varphi^2}{2\pi\xi}\,d\xi.
$$
\end{thm}

The next section is devoted to the proof of this result.

\section{Proof of Theorem \ref{MainResult}}

Let us first give estimates of the trace on $\Gamma\de$ for functions of $W^1(\Omega'\de)$
or $W^{1,p}(\Omega'\de)$, $\frac 32<p<2$.

\begin{lem}
There is a constant $C$, independent of $\eps$, such that~:
\begin{align}
&\|\psi\|_{L^2(\Gamma\de)}\le C\eps^{\frac 12}\,|\ln\eps|^{\frac 12}\,\|\psi\|_{W^1(\Omega'\de)}
\qquad\text{for all }\psi\in W^1(\Omega'\de),\label{Trace1}\\
&\|\psi\|_{L^2(\Gamma\de)}\le C\left(\eps^{\frac 12}\,\|\psi\|_{W^{1,p}(\Omega'\de)}
+\eps^{\frac 43-\frac 2p}\,\|\nabla\psi\|_{\LL^p(\Lambda\de^\delta)}\right)\notag\\
&\qquad\qquad\qquad\text{for all }\psi\in W^{1,p}(\Omega'\de)
\text{ with compact support, $\frac 32<p<2$}.\label{Trace2}
\end{align}
\end{lem}

\begin{proof}
Let $\psi\in{\mathcal C^1}(\overline\Omega'\de)$ with compact support and let
$\h\psi:\h\Omega\to\mathbb R$ defined by
$$
\h\psi(\h{\bs x}) = \psi(\bs F\de(\h{\bs x})),\qquad\h{\bs x}\in\h\Omega.
$$
Let us first prove \eqref{Trace1}. We have
$$
\h\psi(s,0,\theta) = \h\psi(s,1,\theta) - \int_0^1\pd{\h\psi}\xi(s,\xi,\theta)\,d\xi,
\qquad (s,\theta)\in\h\Gamma.
$$
Consequently,
\begin{equation}
|\h\psi(s,0,\theta)|^2\le 2|\h\psi(s,1,\theta)|^2 + 2\left(\int_0^1\pd{\h\psi}\xi(s,\xi,\theta)\,d\xi\right)^2,
\label{Id2}
\end{equation}
and, using the Cauchy--Schwarz inequality and \eqref{estim-param}~:
\begin{align}
|\h\psi(s,0,\theta)|^2 &\le 2|\h\psi(s,1,\theta)|^2
+ 2\left(\int_0^1\frac 1{a\de r\de}\,d\xi\right)\,
\left(\int_0^1a\de r\de\Big|\pd{\h\psi}\xi\Big|^2\,d\xi\right)\notag\\
&\le 2|\h\psi(s,1,\theta)|^2 + 2C_1\,\left(\int_0^1\frac 1{r\de}\,d\xi\right)\,
\left(\int_0^1 a\de r\de\Big|\pd{\h\psi}\xi\Big|^2\,d\xi\right)\notag\\
&\le 2|\h\psi(s,1,\theta)|^2 + C_2\,|\ln\eps|\,\int_0^1 a\de r\de
\Big|\pd{\h\psi}\xi\Big|^2\,d\xi,\label{Id1}
\end{align}
for $(s,\theta)\in\h\Gamma$. Since by \eqref{Int-Gamma},
\begin{align}
&\|\psi\|^2_{L^2(\Gamma\de)} = \eps\int_{\h\Gamma}\alpha\de(s,\theta)|\h\psi(s,0,\theta)|^2\,ds\,d\theta,
\label{Int-Geps}\\
&\|\psi\|^2_{L^2(\Gamma_\delta)} = \delta\int_{\h\Gamma}\alpha_\delta(s,\theta)|\h\psi(s,1,\theta)|^2\,ds\,d\theta,
\label{Int-Gdelta}
\end{align}
with
$$
\alpha_\eps(s,\theta) = a\de(s,0,\theta),\quad \alpha_\delta(s,\theta) = a_\delta(s,1,\theta).
$$
We deduce from \eqref{Id1}, after multiplication by $\eps\alpha\de$ and integration in $s,\theta$,
$$
\|\psi\|^2_{L^2(\Gamma\de)}\le 2\eps\int_{\h\Gamma}\alpha\de|\h\psi(s,1,\theta)|^2\,ds\,d\theta
+ C_2\eps |\ln\eps|\int_{\h\Omega}\alpha\de a\de r\de \Big|\pd{\h\psi}\xi\Big|^2\dhx.
$$
Using \eqref{estim-param} and the estimates $0<C'_3\le\alpha\de$, $\alpha_\delta\le C'_4$, we get
$$
\|\psi\|^2_{L^2(\Gamma\de)}\le C_3\eps\,\delta\int_{\h\Gamma}\alpha_\delta
|\h\psi(s,1,\theta)|^2\,ds\,d\theta
+ C_4\eps\,|\ln\eps|\int_{\h\Omega}a\de r\de\Big|\pd{\h\psi}\xi\Big|^2\dhx.
$$

But \eqref{VariableChange} yields
$$
\|\nabla\psi\|^2_{\LL^2(\Lambda\de^\delta)} = (\delta-\eps)\int_{\h\Omega}
\bigg(\dfrac {r\de}{a\de}\Big(\pd{\h\psi}s-\tau\pd{\h\psi}\theta\Big)^2 + \frac{r\de a\de}
{(\delta-\eps)^2}\Big(\pd{\h\psi}\xi\Big)^2
+ \frac{a\de}{r\de}\Big(\pd{\h\psi}\theta\Big)^2\bigg)\dhx.
$$
Therefore
$$
\|\psi\|^2_{L^2(\Gamma\de)}\le C_3\eps\|\psi\|^2_{L^2(\Gamma_\delta)}
+ C_5\eps|\ln\eps|\,\|\nabla\psi\|^2_{\LL^2(\Lambda\de^\delta)}.
$$
Using the trace inequality and the fact that the support of $\psi$ is compact, we obtain
\begin{align}
\|\psi\|^2_{L^2(\Gamma\de)} &\le (C_3\,C_6\,\eps + C_5\eps|\ln\eps|)\,
\|\nabla\psi\|^2_{\LL^2(\Omega'\de)}\notag\\
&\le C_7\,\eps\,|\ln\eps|\,\|\nabla\psi\|^2_{W^1(\Omega'\de)}.\notag
\end{align}
By density, \eqref{Trace1} follows.

Let us now prove \eqref{Trace2}. We have
\begin{align}
|\h\psi(s,0,\theta)|^2 &= |\h\psi(s,1,\theta)|^2 - \int_0^1\pd{}{\xi}(\h\psi)^2\,d\xi\notag\\
&= |\h\psi(s,1,\theta)|^2 - 2\int_0^1\h\psi\pd{\h\psi}{\xi}\,d\xi.\notag
\end{align}
Multiplying by $\eps$ and integrating in $s,\theta$, we get
$$
\eps\int_0^1\int_0^{2\pi}|\h\psi(s,0,\theta)|^2\,d\theta\,ds
= \eps\int_0^1\int_0^{2\pi}|\h\psi(s,1,\theta)|^2\,d\theta\,ds -
2\eps\int_{\h\Omega}\h\psi\pd{\h\psi}{\xi}\dhx.
$$
Using \eqref{Int-Geps}, \eqref{Int-Gdelta} and \eqref{estim-param}, we get
\begin{equation}
\|\psi\|^2_{L^2(\Gamma\de)} \le C_8\eps\,\|\psi\|^2_{L^2(\Gamma_\delta)} +
C_9\eps\,\left|\int_{\h\Omega}\h\psi\pd{\h\psi}{\xi}\dhx\right|.\label{Id3}
\end{equation}
To estimate the integral in the previous relationship we use the H\"older inequality
$$
\left|\int_{\h\Omega}\h\psi\pd{\h\psi}{\xi}\dhx\right| \le
\Big(\int_{\h\Omega}r\de |\h\psi|^q\dhx\Big)^{\frac 1q}
\Big(\int_{\h\Omega}r\de \Big|\pd{\h\psi}\xi\Big|^p\dhx\Big)^{\frac 1p}
\Big(\int_{\h\Omega}r\de^{1-m}\dhx\Big)^{\frac 1m},
$$
where $q = \frac{3p}{3-p}$ and $m$ is such that $\frac 1p+\frac 1q + \frac 1m=1$,
{\em i.e.,} $m=\frac {3p}{4p-6}$. Using \eqref{dv-db}--\eqref{dv-dt}, we have
\begin{align}
\|\nabla\psi\|_{\LL^p(\Lambda\de^\delta)} = \Bigg(\int_{\h\Omega}
(\delta-\eps)a\de r\de\Bigg(&\frac 1{a\de^2}\Big(\pd{\h\psi}s-\tau\pd{\h\psi}\theta\Big)^2+
\frac 1{r\de^2}\Big(\pd{\h\psi}\theta\Big)^2\notag\\
&+\frac 1{(\delta-\eps)^2}\Big(\pd{\h\psi}\xi\Big)^2\Bigg)^{\frac p2}\dhx\Bigg)^{\frac 1p}.\notag
\end{align}
Using \eqref{estim-param}, we then have
\begin{align}
\left|\int_{\h\Omega}\h\psi\pd{\h\psi}{\xi}\dhx\right| &\le
C_{10}\,\|\psi\|_{L^q(\Lambda\de^\delta)}\,\|\nabla\psi\|_{\LL^p(\Lambda\de^\delta)}\,
\Big(\int_{\h\Omega}r\de^{1-m}\dhx\Big)^{\frac 1m}\notag\\
&\le C_{11}\eps^{\frac{2-m}m}\,\|\psi\|_{L^q(\Lambda\de^\delta)}\,\|\nabla\psi\|_{\LL^p(\Lambda\de^\delta)}.\notag
\end{align}
We note here that $m>2$. Then the imbedding of $W^{1,p}(\Lambda\de^\delta)$ into $L^q(\Lambda\de^\delta)$ implies
$$
\left|\int_{\h\Omega}\h\psi\pd{\h\psi}{\xi}\dhx\right| \le
C_{12}\,\eps^{\frac{2-m}m}\,\|\nabla\psi\|^2_{\LL^p(\Lambda\de^\delta)} =
C_{12}\,\eps^{\frac 53 -\frac 4p}\,\|\nabla\psi\|^2_{\LL^p(\Lambda\de^\delta)}.
$$
Putting this estimate into \eqref{Id3} yields
$$
\|\psi\|^2_{L^2(\Gamma\de)} \le C_8\,\eps\,\|\psi\|^2_{L^2(\Gamma_\delta)} +
C_9\,C_{12}\,\eps^{\frac 83 -\frac 4p}\,\|\nabla\psi\|^2_{\LL^p(\Lambda\de^\delta)}.
$$
Using the trace inequality
$$
\|\psi\|_{L^2(\Gamma_\delta)} \le C_{13}\,\|\psi\|_{W^{1,p}(\Omega'_\delta)},
$$
we get
\begin{align}
\|\psi\|^2_{L^2(\Gamma\de)} &\le C\left(\eps\,\|\psi\|^2_{W^{1,p}(\Omega'_\delta)} +
\eps^{\frac 83 -\frac 4p}\,\|\nabla\psi\|^2_{\LL^p(\Lambda\de^\delta)}\right)\notag\\
&\le C\left(\eps\,\|\psi\|^2_{W^{1,p}(\Omega'\de)} +
\eps^{\frac 83 -\frac 4p}\,\|\nabla\psi\|^2_{\LL^p(\Lambda\de^\delta)}\right)\notag
\end{align}

The conclusion of the lemma follows by density.
\end{proof}

\subsection{Proof of Estimate \eqref{err-u}}

Let $\widetilde w_2\ue=w_2\ue-w_2$. Clearly $\widetilde w_2\ue=u\ue-u$,
$\widetilde w_2\ue\in W^1(\Omega'\de)$ and it satisfies
\begin{equation}
\left\{\begin{aligned}{}
&\Delta\widetilde w_2\ue = 0 &&\qquad\text{in }\Omega'\de,\\
&\pd{\widetilde w_2\ue}{n\de} = -\pd{w_1}{n\de} - \pd{w_2}{n\de}
&&\qquad\text{on }\Gamma\de,\\
&\widetilde w_2\ue(\bs x) = O(|\bs x|^{-1}), &&\qquad |\bs x|\to +\infty.
\end{aligned}\right.\label{Pb-1}
\end{equation}
Using the variational formulation associated with \eqref{Pb-1}, Cauchy--Schwarz inequality and
Estimate \eqref{Trace1}, we deduce
\begin{align}
\int_{\Omega'\de}|\nabla\widetilde w_2\ue|^2\dx
&= \int_{\Gamma\de}\left(\pd{w_1}{n\de} + \pd{w_2}{n\de}\right)\widetilde w_2\ue\,d\sigma\notag\\
&\le \left\|\pd{w_1}{n\de} + \pd{w_2}{n\de}\right\|_{L^2(\Gamma\de)}\|\widetilde w_2\ue\|_{L^2(\Gamma\de)}\notag\\
&\le C\,\eps^{\frac 12}\,|\ln\eps|^{\frac 12}\,\left(\left\|\pd{w_1}{n\de}\right\|_{L^2(\Gamma\de)}
 + \left\|\pd{w_2}{n\de}\right\|_{L^2(\Gamma\de)}\right)
\|\nabla\widetilde w_2\ue\|_{\LL^2(\Omega'\de)}.\label{estim-w2t}
\end{align}
Using \eqref{Estim-dw1dn-inf}, we have
\begin{equation}
\left\|\pd{w_1}{n\de}\right\|_{L^2(\Gamma\de)}\le C\,(\text{meas}\,\Gamma\de)^{\frac 12}\le
C_1\eps^{\frac 12}.\label{estim-dw1dn}
\end{equation}
To estimate $\pd{w_2}{n\de}$, we use standard regularity results for elliptic problems,
see \cite{Dautray-Lions-1}, p. 343, to deduce, since $f\in L^p(\mathbb R^3)$
for $p<2$, that $w_2\in W^{2,p}_{\text{loc}}(\mathbb R^3)$. Then we apply Estimate
\eqref{Trace2} to the function $u=\pd{w_2}{x_i}$, $1\le i\le 3$ with $p=2-\eta$,
$0<\eta<\frac 12$,
$$
\left\|\pd{w_2}{x_i}\right\|_{L^2(\Gamma\de)}\le C\,\left(\eps^{\frac 12}\left\|\pd{w_2}{x_i}
\right\|_{W^{1,p}(\Omega'\de)}
+\eps^{\frac 13-\frac\eta{2-\eta}}\left\|\pd{}{x_i}\nabla w_2\right\|_{\LL^p(\Lambda\de^\delta)}\right).
$$
Since both norms on the right--hand side of the above inequality are uniformly bounded and
since the outward unit normal $n\de$ is uniformly bounded we obtain
\begin{equation}
\left\|\pd{w_2}{n\de}\right\|_{L^2(\Gamma\de)} \le C\eps^{\frac 13-\frac\eta{2-\eta}}.\label{estim-dw2dn}
\end{equation}
Reporting \eqref{estim-dw1dn} and \eqref{estim-dw2dn} into \eqref{estim-w2t} and using the
inequality $|\ln\eps|\le C\eps^{-2\eta}$, we get
$$
\int_{\Omega'\de}|\nabla \widetilde w_2\ue|^2\dx \le C_1\,\eps^{\frac 56-\frac\eta{2-\eta}-\eta}
\,\|\nabla\widetilde w_2\ue\|_{\LL^2(\Omega'\de)}.
$$
Therefore
$$
\|\nabla\widetilde w_2\ue\|_{\LL^2(\Omega'\de)}\le C_2\eps^{\frac 56-\eta}
\qquad\text{for all } \eta>0.
$$
\qed

\subsection{Proof of Estimate \eqref{expr-Le}}

To prove \eqref{expr-Le} we need the following lemmas.

\begin{lem}\label{lem-Le1}
We have for all $\eta>0$,
\begin{equation}
L\ue = \int_{\Omega'\de}|\nabla v|^2\,d\bs x - \int_{\mathbb R^3}fw\dx +
\int_\Sigma(1-\varphi)\left(\pd wn + 2\pd vn\right)\,d\sigma + O(\eps^{\frac 56-\eta}),
\label{lem-Le}
\end{equation}
where $w=w_1+w_2$.
\end{lem}

\begin{proof}
Using the decomposition $u\ue=v+w\ue=v+w_1+w_2\ue$ it follows~:
$$
L\ue = \int_{\Omega'\de\setminus\Sigma}|\nabla v|^2\dx + \int_{\Omega'\de\setminus\Sigma}|\nabla w\ue|^2\dx
+ 2\int_{\Omega'\de\setminus\Sigma}\nabla v\cdot\nabla w\ue\dx.
$$
The estimation of the last two integrals can be achieved as follows. We use \eqref{Pb-w} and the
Green's formula to obtain
\begin{align}
\int_{\Omega'\de\setminus\Sigma}|\nabla w\ue|^2\dx &= -\int_{\Omega'\de\setminus\Sigma}
w\ue\Delta w\ue\dx - \int_{\Gamma\de}w\ue\pd{w\ue}{n\de}\,d\sigma + \int_\Sigma(1-\varphi)
\pd{w\ue}{n\de}\,d\sigma\notag\\
&= \int_{\Omega'\de}fw\ue\dx + \int_\Sigma (1-\varphi)\pd{w\ue}n\,d\sigma.\notag
\end{align}
Similarly, we use \eqref{Singular-Pb} to get
\begin{align}
\int_{\Omega'\de\setminus\Sigma}\nabla v\cdot\nabla w\ue\dx
&= -\int_{\Omega'\de\setminus\Sigma} w\ue\Delta v\dx - \int_{\Gamma\de}w\ue\pd v{n\de}\,d\sigma +
\int_\Sigma(1-\varphi)\pd vn\,d\sigma\notag\\
&= -\int_{\Omega'\de}fw\ue\dx + \int_\Sigma (1-\varphi)\pd vn\,d\sigma.\notag
\end{align}
Then
\begin{equation}
L\ue = \int_{\Omega'\de\setminus\Sigma}|\nabla v|^2\dx - \int_{\Omega'\de}fw\ue\dx +
\int_\Sigma (1-\varphi)\left(\pd{w\ue}n+2\pd vn\right)d\sigma.\label{Le-1}
\end{equation}
We can now estimate the error between the above expression of $L\ue$ and the desired one.
We have, with $w=w_1+w_2$,
\begin{align}
\left|\int_{\mathbb R^3}fw\dx-\int_{\Omega'\de}fw\ue\dx\right| &=
\left|\int_{\mathbb R^3}fw_1\dx - \int_{\Omega'\de}fw_1\dx + \int_{\mathbb R^3}fw_2\dx
-\int_{\Omega'\de}fw_2\ue\dx\right| \notag\\
&\le \left|\int_{\Omega\de}fw_1\dx\right| + \left|\int_{\mathbb R^3}fw_2\dx-\int_{\Omega'\de}fw_2\dx\right|\notag\\
&\qquad\qquad\qquad\quad\! +\left|\int_{\Omega'\de}f(w_2-w_2\ue)\dx\right|\notag\\
&\le \left|\int_{\Omega\de}fw_1\dx\right| + \left|\int_{\Omega\de}fw_2\dx\right|
+ \left|\int_{\Omega'\de}f(w_2-w_2\ue)\dx\right|.\notag
\end{align}
For $1\le p<2$ and $q$ such that $\frac 1p+\frac 1q=1$, we have thanks to Lemma \ref{f}
and since $w_2\in W^{2,p}(\Omega_\delta)\subset L^\infty(\Omega_\delta)$,
\begin{align}
\left|\int_{\Omega\de}fw_2\dx\right| &\le \|f\|_{L^p(\Omega\de)}\,\|w_2\|_{L^q(\Omega\de)}\notag\\
&\le  \|f\|_{L^p(\Omega\de)}\,\|w_2\|_{L^\infty(\Omega\de)}\,(\text{meas}\,\Omega\de)^{\frac 1q}\notag\\
&\le C\eps^{\frac 2q}.\notag
\end{align}
We also have, since $1-\varphi=0$ in a neighborhood of $\partial\Sigma$ and then $w_1 \in H^2(\Omega_{\frac\delta 2})
\subset L^\infty(\Omega_{\frac\delta 2})$,
$$
\left|\int_{\Omega\de}fw_1\dx\right| \le \|f\|_{L^p(\Omega\de)}\,\|w_1\|_{L^\infty(\Omega\de)}
\,(\text{meas}\,\Omega\de)^{\frac 1q}\le C\,\eps^{\frac 2q}.
$$
In addition, since $w_2-w_2\ue=u-u\ue$,
$$
\left|\int_{\Omega'\de}f(w_2-w_2\ue)\dx\right|\le \|f\|_{L^p(\Omega'\de)}\,\|u-u\ue\|_{L^q(\Omega'\de)}.
$$
Choosing $p$ so that $q<\frac{12}5$ and using \eqref{err-u}, we obtain
$$
\left|\int_{\Omega\de}fw_1\dx\right| + \left|\int_{\Omega\de}fw_2\dx\right|
+ \left|\int_{\Omega'\de}f(w_2-w_2\ue)\dx\right|\le C\,\eps^{\frac 56-\eta},
$$
for any $\eta>0$. Now we have to estimate the difference of the two integrals over
$\Sigma$ in \eqref{Le-1} and in \eqref{lem-Le}. From \eqref{Pb-1}, \eqref{eq-w2}
and the identity $\widetilde w_2\ue=w_2\ue-w_2$, we deduce
\begin{align}
\int_\Sigma(1-\varphi)\left(\pd{w\ue}n-\pd wn\right)\,d\sigma &=
\int_\Sigma(1-\varphi)\pd{\widetilde w_2\ue}n\,d\sigma\notag\\
&= \int_{\Omega'\de\setminus\Sigma}\nabla w_2\cdot\nabla\widetilde w_2\ue\dx
+ \int_{\Gamma\de}w_2\left(\pd{w_1}{n\de}+\pd{w_2}{n\de}\right)\,d\sigma.\notag
\end{align}
Then, using estimates \eqref{Trace1}, \eqref{estim-dw1dn} and \eqref{estim-dw2dn}
we get for any $0<\eta\le\frac 12$,
\begin{align}
\left|\int_\Sigma (1-\varphi)\left(\pd{w\ue}n-\pd{w}n\right)\,d\sigma\right|
&\le \|\nabla w_2\|_{\LL^2(\Omega'\de)}\,\|\nabla\widetilde w_2\ue\|_{\LL^2(\Omega'\de)}\notag\\
&\quad+ \|w_2\|_{L^2(\Gamma\de)}\left(\left\|\pd{w_1}{n\de}\right\|_{L^2(\Gamma\de)}
+\left\|\pd{w_2}{n\de}\right\|_{L^2(\Gamma\de)}\right)\notag\\
&\le \|\nabla w_2\|_{\LL^2(\Omega'\de)}\,\|\nabla\widetilde w\ue_2\|_{\LL^2(\Omega'\de)}\notag\\
&\quad + C\eps^{\frac 12}\,|\ln\eps|^{\frac 12}\,\|w_2\|_{W^1(\Omega'\de)}(\eps^{\frac 12}+\eps^{\frac 13-\frac\eta{2-\eta}})\notag\\
&\le C_1\,\left(\|\nabla\widetilde w_2\ue\|_{\LL^2(\Omega'\de)} + C\,|\ln\eps|^{\frac 12}
(\eps+\eps^{\frac 56-\frac\eta{2-\eta}}\right).\label{Int-2}
\end{align}
Using the identity $\widetilde w_2\ue=u\ue-u$ and \eqref{err-u}, we get
\begin{equation}
\left|\int_\Sigma(1-\varphi)\left(\pd{w\ue}n-\pd wn\right)\,d\sigma\right|\le C_2\,\eps^{\frac 56-\eta},
\label{Err-Sigma}
\end{equation}
for any $\eta>0$. Then we obtain the lemma from \eqref{Le-1}--\eqref{Err-Sigma}
\end{proof}

\begin{lem}\label{lem-Le2}
We have
$$
\int_{\Omega'\de\setminus\Sigma}|\nabla v|^2\dx = -\frac{\ell_\gamma}{2\pi}\ln\eps + L' + O(\eps),
$$
where $\ell_\gamma$ is the length of the curve $\gamma$ and
$$
L' = \frac{\ell_\gamma}{2\pi}\ln{\frac\delta 2} +
\frac 1{4\pi^2}\int_{\h\Omega}\left(a_0\xi\theta^2(\h\varphi')^2 +
\frac{\delta^2\xi\tau^2}{a_0}\h\varphi^2\right)\dhx
+ \frac{\ell_\gamma}{2\pi}\int_{\frac 12}^1\frac{\h\varphi^2}\xi\,\dxi.
$$
\end{lem}

\begin{proof}
Using the definition of $v$ and the change of variable $\bs x=\bs F_0(\h{\bs x})$, it follows
$$
\int_{\Omega'\de\setminus\Sigma}|\nabla v|^2\dx = \int_{\Lambda\de^\delta}|\nabla v|^2\dx = A\de^\delta + B\de^\delta
$$
with
\begin{align}
&A\de^\delta = \int_{\h\Omega^\delta\de}\frac{a_0\xi\theta^2}{4\pi^2}(\h\varphi')^2\dhx
+ \int_{\h\Omega^\delta\de}\frac{\delta^2\xi\tau^2}{4\pi^2a_0}\h\varphi^2\dhx,\notag\\
&B\de^\delta = \int_{\h\Omega^\delta\de}\frac{a_0}{4\pi^2\xi}\h\varphi^2\dhx,\notag
\end{align}
where $\h\Omega\de^\delta=(0,1)\times(\frac\eps\delta,1)\times(0,2\pi)$.
Clearly, we can write
\begin{equation}
A\de^\delta = \int_{\h\Omega}\frac{a_0\xi\theta^2}{4\pi^2}(\h\varphi')^2\dhx
+ \int_{\h\Omega}\frac{\delta^2\xi\tau^2}{4\pi^2a_0}\h\varphi^2\dhx + O(\eps).
\label{Ae}
\end{equation}
Since $\h\varphi(\xi)=1$ for $0\le\xi\le\frac 12$, we can write $B\de^\delta$ as
\begin{align}
B\de^\delta &= \int_0^1\!\!\!\int_{\frac\eps\delta}^{\frac 12}\!\!\!\int_0^{2\pi}\frac{a_0}{4\pi^2\xi}\,d\theta
\,d\xi\,ds +
\int_0^1\!\!\!\int_{\frac 12}^1\!\!\!\int_0^{2\pi}\frac{a_0}{4\pi^2\xi}\h\varphi^2\,
d\theta\,d\xi\,ds\notag\\
&=\frac 1{2\pi}\left(\int_0^1|\bs g'(s)|\,ds\right)\int_{\frac\eps\delta}^{\frac 12}
\frac{d\xi}\xi + \int_0^1\!\!\!\int_{\frac 12}^1\!\!\!\int_0^{2\pi}\frac{a_0}{4\pi^2\xi}
\h\varphi^2\,d\theta\,d\xi\,ds\notag\\
&= -\frac{\ell_\gamma}{2\pi}\ln\eps + \frac{\ell_\gamma}{2\pi}\ln\frac\delta 2
+ \int_0^1\!\!\!\int_{\frac 12}^1\!\!\!\int_0^{2\pi}\frac{a_0}{4\pi^2\xi}
\h\varphi^2\,d\theta\,d\xi\,ds\notag\\
&= -\frac{\ell_\gamma}{2\pi}\ln\eps + \frac{\ell_\gamma}{2\pi}\ln\frac\delta 2
+ \int_0^1\!\!\!\int_{\frac 12}^1\!\!\!\int_0^{2\pi}\frac{|\bs g'|}{4\pi^2\xi}
\h\varphi^2\,d\theta\,d\xi\,ds\notag\\
&\qquad
- \int_0^1\!\!\!\int_{\frac 12}^1\!\!\!\int_0^{2\pi}\frac{\delta\xi\kappa\cos\theta}{4\pi^2\xi}
\h\varphi^2\,d\theta\,d\xi\,ds.\notag\\
&= -\frac{\ell_\gamma}{2\pi}\ln\eps + \frac{\ell_\gamma}{2\pi}\ln\frac\delta 2
+ \frac{\ell_\gamma}{2\pi}\int_{\frac 12}^1\frac{\h\varphi^2}\xi\,d\xi.\notag
\end{align}
From this and \eqref{Ae} follows the lemma.
\end{proof}

Estimate \eqref{expr-Le} follows immediately by combining Lemmas \ref{lem-Le1} and
\ref{lem-Le2}.

\bibliographystyle{unsrt}

\end{document}